\input amstex

\documentstyle {amsppt}
\NoRunningHeads

\topmatter

\title
About approximate and dual Lie algebras
\endtitle

\author V.V. ~ Gorbatsevich
\endauthor

\endtopmatter

\document

\head Introduction o
\endhead

In this article the concept of an approximate Lie algebra, depending on small parameter, and which is naturally arising when studying approximate symmetry (about which see, for example, fundamental article [1]) of the differential equations, is considered. At the same time it is natural to pass to more general concept of dual or $D_2$-Lie algebras, for which many results are represented in more convenient form. The main purpose of this article --- to clear some algebraic concepts, which are used when studying approximate solutions of the equations and their approximate symmetry. The author hopes that the results, which are contained in this article, will be useful when studying, in particular, of Lie algebras of infinitesimal approximate symmetry of differential equations.

In Introduction the way, which led to use of the concept of dual or $D_2$-Lie algebras, will be described.

The concept of the approximate Lie algebra, which we will study, was entered for a long time. In the beginning it was used only for Lie algebras of vector fields (or --- it was very frequent --- for the corresponding differential operators) and as a separate object of the research is allocated recently (see, for example, [2,3,4]) and further underwent quite detailed studying. At the same time the basic idea of consideration of functions and vector fields was to approximate them --- up to $o(\epsilon)$, what is equivalent to their linearization. For example, smooth (or even analytical) function $F(x,\epsilon)$ is replaced with the function of the type of $f_0 (x)+\epsilon f_1 (x)$, where $f_0,f_1$ --- some smooth (or analytical) functions. Such replacement will be designated further by symbol $\approx$. Similarly, for smooth (or analytical) vector field $X (x, \epsilon) $, depending on parameter $\epsilon$, we have approximate equality $X(x,\epsilon) \approx X_0(x) + \epsilon X_1 (x)$ , where $X_0 (x), X_1(x)$ --- vector fields. The commutation of such approximate vector fields can be calculated by such formulas:

$$ [X, Y] \approx [X_0+\epsilon X_1, Y_0+\epsilon Y_1] \approx [X_0, Y_0] + \epsilon ([X_0, Y_1]+ [X_1, Y_0]) $$

Here the bilinearity of operation of such commutation and skew-symmetry are remain, Jacobi's identity is carried out too. So we receive some Lie algebra. But if to consider the standard commutation of the linearized vector fields, then the result will not be, generally speaking, linear on $\epsilon$ any more. Therefore the given above commutation is not standard operation of commutation of vector fields --- it is made in a different way. We come to the new concept, in which in the uniform form all Lie algebras, depending on parameter $\epsilon$, are combined. This object is called as approximate Lie algebra and now we will consider several different approaches to consideration of such objects.

Construction 1

Let $\Phi$ be some Lie algebra of vector fields (or their jets) in some neighborhood of the point $x \in \bold R^n$. Let's consider the set of vector fields of the form $X + \epsilon Y$, where $X, Y \in \Phi$, and $\epsilon$  is some parameter($\epsilon^2=0$). This set we will designate by $\Phi (\epsilon)$ or we will present it in the form $\Phi + \epsilon \Phi$. Generally speaking, this set of vector fields (or of jets) does not form a Lie algebra, as operation of commutation can lead to the vector fields which are not lying in $\Phi (\epsilon)$. But if we define on $\Phi (\epsilon)$ a commutation operation in a different way (as it is made above for approximate vector fields), then we will get, as it is easy to be verified, the Lie algebra (but it will not be the a Lie algebra of vector fields any more). Such Lie algebras are called approximate one's in the articles stated above.

Except the approximate Lie algebras, which are already constructed by such a way, it is necessary to consider also some of their subalgebras --- those, which are stable under multiplication on $\epsilon$. Such Lie subalgebras correspond to Lie groups of approximate symmetry of differential equations, which studying initiated introduction of the concept of the approximate Lie algebra.

Let's note also, that consideration of expressions of the similar type $\bar a + \omega \bar a^0$, where $\bar a,  \bar a^0$  are three-dimensional vectors, and $\omega$ be the symbol, for which $\omega ^2 = 0$, is used in mechanics as subject to screw calculation.

Construction 2

Let $\Phi$ be any Lie algebra. Let's consider two-dimensional algebra $D_2$ (algebra of dual numbers of the form $a+\epsilon b$; sometimes it is called Studi's algebra and its elements---Studi's numbers) over $\bold R$ with generators 1 and $\epsilon$, where $\epsilon ^2 = 0$. Let's note that this algebra $D_2$ (which is commutative, associative and with a unit) has the exact matrix representation, which matrices are $\left( \smallmatrix a&b \\0&a \endsmallmatrix \right )$, where $a, b \in \bold R$. It is also interesting to note, that unlike the field $\bold R$. which group of automorphisms is trivial, and field of complex numbers  $\bold C$, which group of automorphisms is isomorphic $\bold Z_2$ (with complex conjugation as the only nontrivial automorphism), the algebra $D_2$ has one-dimensional group of automorphisms---it consists of transformations of the form $a+\epsilon b \to a+\epsilon \alpha b$ for any nonzero $\alpha \in \bold R$.

Let's denote $\tilde{\Phi } = \Phi \otimes D_2$ - it is the tensor product of algebras $\Phi$ and $D_2$. It is possible to write down elements of $\tilde{\Phi}$ in the form of $a+\epsilon b$, where $a,b \in \Phi$. Operation of commutation in the algebra $\tilde{\Phi}$ is set naturally (with commutation on the first and multiplication on the second tensor factors), so we get on $\tilde{\Phi } $ the structure of the Lie algebra. Up to isomorphism of Lie algebras $\tilde{\Phi }$ can be considered as the semidirect sum $\Phi +_ {ad} \Phi$ of the subalgebra, which is isomorphic to $\Phi$, and an Abelian ideal, corresponding to the adjoint action $ad$ of Lie algebra $\Phi$ on itself as on the vector space. It is clear, that if $\Phi$ be a Lie algebra of vector fields, then by such construction we get the approximate Lie algebra from the Construction 1.

This Construction 2 can be generalized. For example, we may consider the algebra of the form $D_p = \bold R [\epsilon] / <\epsilon ^p> $, where $\bold R (\epsilon) $ is the algebra of polynomials and $ <\epsilon ^p> $ is the ideal, generated by element $\epsilon^p$. We get an associative and commutative algebra of dimension $p$ with unit, sometimes it is called as algebra of plural numbers. For $p=2$ this algebra was already entered above. Let's consider the tensor product $\Phi \otimes D_p$ with natural operation of commutation--- we  get here the Lie algebra, which can be considered as the approximate Lie algebra when using approximations of order $p$ by means of $o (\epsilon^p) $. It is presented in the form of the semidirect sum of subalgebra, isomorphic to $\Phi$, and the nilpotent ideal of the form $\epsilon \Phi + \epsilon^2 \cdot \Phi \dots + \epsilon^{p-1} \cdot \Phi$.

The Lie algebras, constructed above by us , are graduated. It is possible to study for the bigger expansion of the concept of the approximate Lie algebra also more general graduated Lie algebras---for example, any 2-graduated Lie algebras of the type of $L = L_0 +L_1$ (which can be written down, by analogy with the aforesaid, in the form of $L_0 + \epsilon L_1)$, where $L_0$ --- some Lie subalgebra, and $L_1$ is the Abelian ideal in $L$. It is possible to consider also some longer graduation (they arise above for $p> 2$).

After the description of two formally different (but very close in principle) approaches to the concept of the approximate Lie algebra, we will pass directly to the statement of our approach to this concept.

By approximate Lie algebras (in more general case, and not only for Lie algebras, connected with vector fields) it would be possible to call, first of all, Lie algebras of the form $L \otimes D_2$ --- tensor products of ordinary Lie algebras $L$ and two-dimensional algebra $D_2 = <1, \epsilon> $ (the algebra of dual numbers). Let's note that such Lie algebras can be considered as $D_2$-modules. Moreover, approximate Lie algebras of symmetry of approximate differential equations always have a structure of $D_2$-modules too. It follows from such simple fact: if the vector field $X$ is an infinitesimal symmetry of a differential equation, that and vector field $\epsilon X$--- too. Thereby the Lie algebra of infinitesimal symmetry is invariant under multiplication on elements of algebra $D_2$.

Except Lie algebras of the type of $L \otimes D_2$, it is possible (and it is necessary) to consider --- for the purpose of their further applications --- and their subalgebras too. Generally speaking, any Lie subalgebra in $L \otimes D_2$ is not closed  under the multiplication on $\epsilon$ (and therefore is not a $D_2$-Lie subalgebra ). For example, for any element $X \in L$ spanned  on it one-dimensional abelian Lie subalgebra in $ L \otimes D_2$ do not $D_2$-invariant. But if we are interested in Lie algebras of approximate symmetries, then it is natural to consider only such Lie subalgebras, which are $D_2$-modules. Also it is possible to call them as approximate Lie algebras (in that sense, what arises when we study approximate differential equations and their symmetries). Thereby we come to need to consider some class of the Lie algebras, defined over $D_2$. Namely, as the next approach to our final definition of approximate Lie algebras,  it make sense to call as approximate subalgebras any invariant under multiplication on $\epsilon$ Lie subalgebras in Lie algebras of the type of $L \otimes D_2$ for various Lie algebras $L$ (in the theory of approximate symmetry of differential equations it is used as $L$ the Lie algebras of smooth or even analytical vector fields). But to separate this general concept from that, which is connected with vector fields, we will call such Lie algebras--- which are  $D_2$-modules - as $D_2$-Lie algebras (also the name "dual Lie algebras" is admissible, though it is a little ambiguous).

So, our main object, which we will also study further, is any $D_2$-Lie algebra. The linear operator of multiplication on $\epsilon$ for such Lie algebras we will denote by $\Cal E$ (actually it determines $D_2$-structure). At the same time, it is natural to study also the corresponding Lie subalgebras, ideals, homomorphisms, etc. Let's note that as $D_2$-module the approximate Lie algebra can be trivial (i.e. operator $\Cal E$ for it can be zero). Thereby ordinary Lie algebras, for which $\Cal E=0 $, can be considered as $D_2$-Lie algebras. However, further we will be interested generally in cases when $\Cal E \ne 0$. Moreover, sometimes it is useful to consider only "nondegenerate" $D_2$-structures, when the rank of operator $\Cal E$ is equal to the half of the dimension of the Lie algebra (which, at the same time, is supposed, naturally, to be even). It is clear, that it will be in only case when the Lie algebra $L$ is a free $D_2$-module. However the class of such nondegenerate Lie algebras is not closed concerning transition to Lie subalgebras (or even only to ideals), or under transitions to quotient algebras by ideals. Generally, upon such transitions the linear operator $\Cal E$ can be "degenerate", in particular---be zero (when we have the trivial $D_2$-module and the corresponding  $D_2$-Lie algebra is the ordinary real Lie algebra without additional nontrivial $D_2$-structure). Further in article some main properties of any $D_2$-Lie algebras will be studied in detail. In particular, it gives a lot of information about approximate Lie algebras, that can be useful for researches in which these concepts arose.

As algebra $D_2$ is not a field (the division by any numbers of type $\epsilon b$ is impossible), the normal concept of the basis for its modules is inapplicable. By essential basis (or $D_2$-basis) for $D_2$-Lie subalgebra we understand the minimal set of the vectors, generating this subalgebra as $D_2$-module. Speaking in other words, it is such minimal set of vectors, which together with their multiplications on $\epsilon$ generates the Lie algebra as the vector space over $\bold R$. Such bases were entered in ([5], p. 41; there they are called essential parameters). Let's note that only for Lie algebras over $D_2$ (and for $D_2$-Lie subalgebras) the concept of essential basis makes sense. The number of elements in it (the number of essential parameters) is some analog of the dimension for the considered Lie algebra $L$ over $D_2$ (here it is possible to use the term "$D_2$-dimension" and to designate this number as $d(L)$). For example, $D_2$-Lie subalgebra in the Lie algebra of the type $L \otimes D_2$, generated by the element of type $\epsilon X$ (where $ X \in L$) is one-dimensional over $\bold R$, but Lie subalgebra, generated by the element $X$, is two-dimensional (here $X=X+\epsilon \cdot 0 $).

If the vector space $V$ is some finite-dimensional free $D_2$-the module, then $V$ is isomorphic to $D_2^m$ for some natural $m$. In this case it is natural to consider in $V$ $D_2$-bases (systems, free generating the module) and linear operators are defined by matrices (which elements belong to $D_2$).

Let's note also, that the concept of the characteristic polynomial is badly adapted for use in $D_2$-the situation. First, the notion of the basis is absent (in normal sense of this word), and therefore it is difficult to use the concept of the matrix. As for polynomials and their roots, then in general there are a lot of surprises. For example, the polynomial can have infinitely many roots---for example, for the polynomial $ (z-a) ^2$ for the fixed $a \in \bold R$ any dual number of the form $a+\epsilon b$ with an arbitrary value of the real number $b$, will be the root of this polynomial.

Here the Introduction, having generally methodological focus, comes to an end and we pass to the studying of the object, entered above--- the category of $D_2$-Lie algebras. We will consider generally real Lie algebras, though absolutely similar considerations can be carried out for complex $D_2$-Lie algebras (note that algebra $D_2$ can be considered also over $\bold R$ and over $\bold C$).

\head \S 1 Some properties of $D_2$-Lie algebras
\endhead

Here some properties of finite-dimensional $D_2$-Lie algebras (i.e. such real Lie algebras, which are $D_2$-modules) will be considered. Let's remind that $D_2$ is the algebra of dual numbers (which can be considered as result of Cayley-Dixon doubling for the algebra of real numbers) and therefore it is possible to call the corresponding Lie algebras as dual Lie algebras. It is possible to consider also algebra of dual numbers over $\bold C$ (several times it will be stated below). We will also be limited here to consideration of finite-dimensional Lie algebras (though some statements, given below, are true also without this restriction).

On any real Lie algebra $L$ the structure of $D_2$-Lie algebra is set by the linear operator $\Cal E$ (which we will often consider as the operator of multiplication by some special element $\epsilon \in D_2)$, for which  there are two following properties:

1. $ [\Cal E X, Y]= [X, \Cal E Y]=\Cal E [X, Y] $

On another way this condition can be expressed so: the operator of the $D_2$-structure commutes  with operators of the adjoint representation of the Lie algebra $L$.

2. $\Cal E^2 =0$

In suitable (Jordan) basis the matrix of the linear operator $\Cal E$ is the direct sum $\oplus J_2 (0) \oplus 0$ of several Jordan cells $J_2(0)$ or order 2, corresponding to the zero eigenvalue, and some zero square matrix (designated above as $0$; it can be absent in this decomposition). The number of Jordan cells $J_2(0)$ in this decomposition is equal to the rank of the matrix $\Cal  E$, and also to the dimension of the image of this linear operator.

Trivial examples of $D_2$-(or dual) Lie algebras are ordinary real Lie algebras, for which the corresponding linear operator $\Cal E$ equals to zero.

The simplest nonstandard examples of $D_2$-Lie algebras are Lie algebras of the type of $L \otimes D_2$--- tensor products of any real Lie algebra $L$ and the algebra $D_2$ (commutation operation is defined naturally here: $ [X \otimes \alpha, Y \otimes \beta] = [X, Y] \otimes (\alpha \beta) $). $D_2$-Lie algebras, which we get by such construction, are even dimensional, if to consider them as real Lie algebras. However there are also odd-dimensional $D_2$-Lie algebras---for example, such is any one-dimensional Lie subalgebra in any $D_2$-Lie algebra, spanned by the element $\epsilon X$. The specified Lie algebras $L \otimes D_2$ are, obviously, free $D_2$-modules. The corresponding linear operator $\Cal E$ has the greatest possible rank (which equals to dimension of the initial Lie algebra $L$).

Odd-dimensional Lie algebras can not be obtained by means of the specified operation of the tensor product (which it is possible to call as operation of the dualization of the initial Lie algebra). However there are also even dimensional $D_2$ - Lie algebras which are not result of the dualization of any Lie algebra. For example, such is the Abelian Lie algebra $\bold R^2$ (with trivial $D_2$-structure). But there are such examples and in the case of nontrivial $D_2$-structures. Here we have an  elementary (having the minimal possible dimension) example - it is the four-dimensional Lie algebra $L_4$, having basis $X, Y, \epsilon X, \epsilon Y$ and the defining bracket $[X, Y]=\epsilon X$ (other nontrivial brackets are zero). It is easy to understand that, as the real Lie algebra, this $L_4$ is isomorphic to the Lie algebra $n(3, \bold R)  \oplus \bold R$, where $n (3, \bold R)$ is the three-dimensional nilpotent Lie algebra (which in suitable basis $U, V, W$ is defined by the bracket $[U, V] =W$). If there is such Lie algebra $P$, that $L$ is isomorphic to $P \otimes D_2$, then the Lie algebra $P$ would have to be two-dimensional. But as is well-known, that there are only two---up to isomorphism, as it is natural---real two-dimensional Lie algebras: they are one Abelian and one solvable (which is not nilpotent). So the Lie algebra $P \otimes D_2$ has to be or Abelian or not nilpotent. But ours Lie algebra $L$ is nilpotent and therefore the assumption of existence of the Lie algebra $P$ leads to the contradiction.
The same conclusion can be obtained more simply (and not based on classification). The matter is that $D_2$-structure on the considered Lie algebra is degenerated (the rank of operator $\Cal E$ is equal to 1). But Lie algebras of the type $L \otimes D_2 $ always have due to their  construction a nondegenerate $D_2$ - structure. Therefore the considered Lie algebra $L_4$ is not the dualization of any Lie algebra. If for $D_2$-Lie algebra the structure is nondegenerate, then examples of this sort can be constructed too, but they will be more bulky (some arguments  in this directiont will be given below).

For any $D_2$-Lie algebra $L$ by $L_{\bold R}$ we will denote the realification of the Lie algebra $L$, i.e. it is Lie algebra $L$ considered over field $\bold R$ (ignoring its $D_2$ - structure). As it will be shown below, many properties of $D_2$-Lie algebras it is convenient to prove, using the known results, applied to real Lie algebras $L_{\bold R }$. However there are also results specific to $D_2$-Lie algebras.

Let $L$ be some finite-dimensional $D_2$-Lie algebra. Let's consider the linear subspace $P = Im (\Cal E) $---the image of the linear operator $\Cal E$, and by $Q$ we will denote its kernel $Ker (\Cal E) $. The condition $\Cal E^2 =0$ is equivalent to the inclusion  $P \subset Q$. It is obvious that the subspace $P$ is zero in only case, when the operator $\Cal E$ is zero. Let's note that number of essential parameters  $d (L)$ for $D_2$-Lie algebra (mentioned above) equals to the real dimension of the vector space $L/P$. We pass to detailed studying of some properties of $D_2$-Lie algebras.

\proclaim {Proposition 1} Let $L$ be some $D_2$-Lie algebra and $U$ be any linear subspace in $L$. Then $\Cal E (U)$ is the Abelian subalgebra in $L$.
\endproclaim

\demo {Proof} We have $ [\Cal E(U), \Cal E (U)] = \Cal E^2 ([U, U]) = \{ 0 \} $.
\enddemo

For any linear subspace $U \subset L$ we will put $\hat U = U+\Cal E U$. Actually $\hat U$ equals to $D_2 \cdot U$ --- it is minimal $D_2$-subspace containing $U$ (or $D_2$- saturation of the subspace $U$).

\proclaim {Proposition 2} Let $U$ be some ideal in $D_2$-Lie algebra.
Then

(i) $\Cal E (U)$ is the ideal in $L$

(ii) $\hat U$ is the ideal in $L$ (and $ [\hat U, L] \subset U$)
\endproclaim

\demo {Proof}
1. We have $ [\Cal E(U), L] = \Cal E ([U, L])$. But $U$ is the ideal in $L$ and therefore $ [U, L] \subset U$, so we get $ [\Cal E(U), L] \subset \Cal E(U)$, i.e. $\Cal E (U)$ is the ideal in $L$.

2. We have $ [\hat U, L] = [U + \Cal E U, L] \subset [U, L] + [\Cal E U, L] $. But $ [U, L] \subset  U$ (because $U$ is the ideal), and also
$ [\Cal E U, L] = [U, \Cal E L] \subset U$. We get $ [\hat U, L] \subset \hat U$, i.e. $\hat U$ is the ideal in $L$. Also $ [\hat U, L] \subset U$ \enddemo

In particular, it follows from Propositions 1 and 2, that subspace $P$, introduced above, is the Abelian ideal in $L$. This statement, as well as the Corollary 1, given below, was noted for the first time in [3]. Let's note that in the specified article some designation is not standard for the theory of Lie algebras --- it is used the notation of decomposition $L = L_0 \oplus L_1$, using designations for the direct sum, though actually the speech there implicitly goes about decomposition to disjoint union of two subsets, from which one---$L_0$---is not the subspace.

\proclaim { Corollary 1 } If $L$ is some $D_2$-Lie algebra and $\Cal E \ne  $0, than $L$ as the Lie algebra over $\bold R$ will never be semi-simple.
\endproclaim

\demo {Proof}. If $\Cal E \ne  0$, then the  subspace $P$  is an nontrivial Abelian ideal, that makes semi-simplicity of the Lie algebra $L$ impossible.
\enddemo

In addition to the statement (ii) of the Proposition 2 we will prove

\proclaim {Proposition 3} Let $U$ be some linear subspace in $D_2$-Lie algebra $L$. Then $\hat U$ is the $D_2$-Lie subalgebra in $L$.
\endproclaim

\demo {Proof} We have $ [\hat U, \hat U] = [U + \Cal E U, U+\Cal E U] \subset [U, U] + [\Cal E U, U] + [\Cal E U, \Cal E U] $. But $ [\Cal E U, U]=\Cal E [U, U] \subset \hat U$ and $ [\Cal E U, \Cal E U]=\Cal E^2 [U,  U]=0$, and therefore $ [\hat U, \hat U] \subset \hat U$ , i.e. $\hat U $ is the Lie subalgebra.
\enddemo

Concepts of the solvable and nilpotent Lie algebras are defined for the class of $D_2$-Lie algebras naturally (by the sequences of the corresponding commutants). It is easy to be convinced that all members of the central series (upper and lower) and the members of commutants series for $D_2$-Lie algebra are $D_2$-Lie subalgebras too. For example, we will show it for the center $Z(L)$: let $X \in Z(L)$, then for any element $Y \in L$ we have $ [\Cal E X, Y] = \Cal E [X, Y] =0$, therefore $\Cal E X \in Z(L)$.

More general statement can be similarly proved

\proclaim {Proposition 4} Let $L$ be some $D_2$-Lie algebra, and $U$ is some $D_2$-Lie subalgebra in it. Then the centralizer $Z_L (U)$ and the normalizer $N_L(U)$ of this Lie subalgebra are $D_2$-Lie subalgebras in $L$.
\endproclaim

Further, the concept of the solvable radical is naturally defined (it is the maximal solvable $D_2$-ideal). As it appears, $D_2$-radical in $D_2$-Lie algebra equals to the ordinary radical (in the Lie algebra $L_{\bold R } $). The situation for the nilradical (the maximal nilpotent ideal of the Lie algebra of $L_{\bold R } $ is similar.

\proclaim { Lemma 1 } Let $L$ be some $D_2$-Lie algebra, and $R$ is the  radical (the maximal solvable ideal) in the Lie algebra $L_{\bold R } $. Then $R$ is $D_2$-radical (the maximal solvable $D_2$-ideal) in the Lie algebra $L$. If $N$ is the nilradical in $L_{\bold R } $, then it is also $D_2$-nilradical.
\endproclaim

\demo {Proof} We will consider Lie algebra $\hat R = R + \Cal E R$. The ideal $R$ will be an $D_2$-ideal in only case, when $\hat R=R$. But $\Cal E (R)$ contains, due to the Proposition 1, in the Abelian ideal $P$ of the Lie algebra $L$. We have $P  \subset R$ and therefore $\hat R \subset R$. Therefore $\hat R = R$.

For the case of the nilradical the proof is much simpler. It, as we know, consists of the radical's elements, for which the operator of the adjoint representation, if restricted on the radical, has only zero eigenvalues. But by multiplication by elements from $D_2$ nonzero eigenvalues cannot appear. Therefore $N=D_2 \cdot N$.
\enddemo

\proclaim { Lemma 2 } Let $L$ be some $D_2$-Lie algebra and $N$ --- its nilradical. Then $N$ contains the Abelian ideal $P = Im \Cal E$. In particular, $\text {dim} N \ge \text {dim} P$.
\endproclaim

\demo {Proof} Due to the Proposition 1, subspace $P$ is the Abelian ideal in $L$ and therefore contains in the radical of the Lie algebra $L$. Let $X \in P$ be any element. Then $ad_X(L) \subset P $ and therefore $ ad^2_X (L) =\{0\}$, i.e. operator $ad_x$ is nilpotent. But therefore $X$ contains in $N$.
\enddemo

The statement of the Lemma 2 can be deduced by a different way: any nilpotent ideal of the Lie algebra (which is contained always in its radical), contains also in its nilradical. It follows from the known description of the nilradical in the solvable Lie algebra, that it is maximal among all nilpotent ideals of the radical.

There are some difficulties with the concept of semisimple $D_2$-Lie algebras and semisimple $D_2$-Lie subalgebras. For them $\Cal E = 0 $ (see the Proposition 1). But further for theory of $D_2$-Lie algebras a lot of things go to as usual---for example, Levi's decomposition $L=S+R$ takes place. At the same time, the radical $R$ will be $D_2$-subalgebra, and the semi-simple part $S$---not always (because Lie subalgebra $\Cal E (S)$---if it is nonzero---never contains in $S$). As the such example we will consider the Lie algebra $S \otimes D_2$, where $S$ --- some semi-simple real Lie algebra. Then follows from the construction of the tensor product, that $\Cal ES \ne \{0\}$, but $\Cal E S \cap S = \{ 0 \}$ (it is the typical case, see the Proposition  5 below).

Let's study in more detail the situation, connected with semi-simple subalgebras in $D_2$-Lie algebra $L$.

\proclaim {Proposition 5} Let $S$ be some semi-simple Lie subalgebra in the finite-dimensional $D_2$-Lie algebra $L$ (more precisely, in $L_{\bold R } $). Then

1. $\Cal E S \cap S = \{ 0 \} $

2. Lie algebra $\hat S$ has Levi's decomposition $S+A$---it is the semidirect sum of the Lie subalgebra $S$ and some Abelian ideal $A \subset \hat S$.
\endproclaim

\demo {Proof}

1. We have $ [\Cal E (S), S] = \Cal E ([S, S]) = \Cal E (S)$ (because
$[S, S] = S$ for any semi-simple Lie algebra $S$). From this it follows, that subspace $\Cal E (S)$ is invariant under action of the Lie algebra $S$, induced by the adjoint representation of the Lie algebra $L$. From another point of view, $\Cal E (S)$ contains in the Abelian ideal $P=\Cal E(L)$ in the Lie algebra $L$ (see the Proposition 1 above). As the semi-simple Lie algebra $S$ has no nontrivial Abelian ideals, we see that $\Cal E S \cap S = \{ 0 \} $, that proves our first statement.

2. Let's consider the Lie subalgebra $\hat L \subset L$. By the definition it is the sum of two subspaces $S$ and $A = \Cal E (S)$ having, as soon as it was proved above, trivial intersection. At the same time the subspace $A$ is the Abelian Lie algebra (as shown above), and it is the Lie subalgebra, which is invariant under the action of $S$ . All this together means that Lie subalgebra $\hat S$ is the semidirect sum of the Lie subalgebra $S$ and the Abelian ideal $A$ (being also its radical).
\enddemo

It is possible to suppose, that if $D_2$- structure on Lie algebra $\hat S$ has the maximum rank, then $\text {dim} A \ge \text {dim}S$. Equality takes place, for example, if $ L=S \otimes D_2$, where $S$---any semi-simple Lie algebra.

Let's consider now Levi's decomposition $L_{\bold R} = S+R$ for realified $D_2$-Lie algebras. Due to the Proposition 5 $\Cal E (S)$ always has the trivial intersection with $S$, i.e. the Lie subalgebra $S$ is far from being invariant under multiplication by elements from $D_2$. This $S$ is isomorphic to quotient $D_2$-Lie algebra (including nondegenerate one's) for the $D_2$-ideal (radical), perhaps degenerate.
For nilpotent $D_2$-Lie algebra its center $Z(L)$ is always nontrivial. Therefore for nilpotent $D_2$-Lie algebras, by using standard induction on dimension, we come to the $D_2$-analog of the Engel theorem:

\proclaim {Proposition 6} Let $L \subset gl(V)$ be some linear $D_2$-Lie algebra, consisting of nilpotent elements (here $V$ is some free $D_2$-module). Then there is such $D_2$-basis in $V$, in which all matrices of the Lie algebra $L$ have the  nilpotent triangular form.
\endproclaim

\demo {Proof} Due to the classical Engel theorem, in $V$ exists vector $X_1$, which is nullified  by all elements from $L$. As $V$ is a free $D_2$-module, the element $X_1$ can be extended to a free basis of $V$. Now we take quotient space $V /  <X> $ and apply induction on dimension.
\enddemo

Further, for solvable $D_2$-Lie algebras some analog of the theorem of Lie takes place (which is formulated usually for complex Lie algebras, but also for real Lie algebras some analog takes place). The classical theorem of Lie for complex linear (or with the fixed linear representation) solvable Lie algebras $L$ claims, that there is one-dimensional invariant subspace in $L$. From this fact by induction it is deduced that all matrices of this linear Lie algebra $L$ in some basis have the triangular form (it   is the statement of Lie's theorem). For real solvable Lie algebras the one-dimensional invariant subspace not always exists, but always there is no more than two-dimensional invariant subspace. Further by induction it is proved, that such Lie algebra has in suitable basis a  quasitriangular form (over $\bold C$ --- triangular) --- on diagonal there are cells of orders 1 or 2, and below---zero elements). It appears that for $D_2$-Lie algebras a similar statement is valid, because it is possible to prove that the specified invariant subspace can be chosen as $D_2$-invariant. For simplicity of the statement we will prove this statement only in case of the adjoint linear representation.

\proclaim {Proposition 7} Let $L$ be a real solvable finite-dimensional $D_2$-Lie algebra. Then in $L$ there is an Abelian $D_2$-ideal of dimension 1 or 2. If $L$ is a complex $D_2$-Lie algebra, then exists $D_2$-ideal of dimension 1.
\endproclaim

\demo {Proof} We will consider the realification $L_{\bold R } $ of the Lie algebra $L$. For the beginning we will assume, that the linear operator $\Cal E$ is nontrivial. Then $P = Im \Cal E$  is the Abelian ideal of positive dimension. Due to the real analog of the theorem of Lie, stated above, in $P$ there is some $L$-invariant subspace of $W$ of dimension (real) 1 or 2. As $P$ is the Abelian ideal, this subspace of $W$ also is the Abelian ideal. At the same time $\Cal E (P) =\{ 0 \}$ (as $\Cal E^2 =0 $). Therefore $\Cal E(W) = \{0\}$ and $W$ is the $D_2$-ideal of dimension (real) 1 or 2.

If operator $\Cal E$ is trivial, then the statement, necessary to us, is  the real analog of the Theorem of Lie, mentioned above.

For the case of the complex Lie algebra the reasoning is similar, only the initial subspace $W$ is one-dimensional.

Thereby the Proposition 7 is completely proved.
\enddemo

Let's show, that for real $D_2$-Lie algebras the one-dimensional ideal not always exists. For this purpose we will consider the Lie algebra $E(2) \otimes D_2$, where $E(2)$ is  the Lie algebra of motions of the Euclidean plane (it is isomorphic to the Lie algebra $so (2)+\bold R^2$). In the Lie algebra $E(2)$ there are no one-dimensional ideals. But then it is easy to check, that there are no one-dimensional ideals also in $D_2$-Lie algebra $E(2) \otimes D_2$.

Let's consider now for $D_2$-Lie algebras the analog of the theorem of Ado about existence of faithful finite-dimensional linear representation for any real finite-dimensional Lie algebra.

Let $V$ be some vector space over $\bold R$. By $V (\epsilon)$ we will denote the vector space $V \otimes D_2$ over algebra $D_2$.
Above such construction was used for introduction of the concept of approximate Lie algebras.

For Lie algebra $gl_n (\bold R)$ of linear transformation of the vector space $\bold R^n$ we will consider corresponding $D_2$-Lie algebra $gl_n (D_2) = gl_n (\bold R) (\epsilon)$ --- it is the Lie algebra of linear transformations of $D_2$-vector space $D_2^n$ (free $D_2$-module). Let's note that it is nondegenerate $D_2$-Lie algebra. It can be presented as consisting of cells of order 2, each of which has the form  $\left ( \smallmatrix x&y \\ 0&x \endsmallmatrix \right)$.
Further, similarly it is possible to introduce other analogs of classical matrix Lie algebras. By  $N(n, D_2)$ we will denote the nilpotent Lie algebra, consisting of all nilpotent upper triangular matrixes with elements from $D_2$ (it is isomorphic to $N (n, \bold R) \otimes D_2$) and it can be written down in the form of the set of block nilpotent matrixes which nonzero blocks are matrices of the form
$\left (\smallmatrix x&y \\ 0&x \endsmallmatrix \right) $. Similarly, by $T (n, D_2) $ we will denote the set of upper triangular matrices with elements from $D_2$---it is solvable (and even triangular) Lie algebra, which matrix representation consists of block matrices of order 2. Let's remind that the real Lie algebra is called triangular (sometimes---completely solvable), if all characteristic roots of operators of its adjoint representation are real.

\proclaim {Proposition 8} Let $L$ be any finite-dimensional Lie algebra over $D_2$. Then it is isomorphic to some $D_2$-Lie subalgebra in $gl_n (D_2)$. In other words, the Lie algebra $L$ has an exact finite-dimensional linear representation over $D_2$.

If $L$ is nilpotent, then it is isomorphic to some $D_2$-Lie subalgebra in $N(n, D_2)$. And when $L$ is triangular (if we consider it as the real Lie algebra), then it is isomorphic to the $D_2$-Lie subalgebra in $T(n,D_2)$.
\endproclaim

\demo {Proof} We will consider $L$ as the Lie algebra $L_{\bold R } $ over $\bold R$ . Due to the classical theorem of Ado there is an embedding  of this Lie algebra $L_{\bold R } $ into Lie algebra $gl (n, \bold R)$ for some $n \in \bold N$. It is clear, that there is the natural embedding $L \subset L_{\bold R } (\epsilon) $. But then we get the $D_2$-Lie algebra $L$ embedding into $gl(n, \bold R) (\epsilon) $.

For cases, when $L$ is nilpotent or triangular, the corresponding statements from the Proposition 8 are proved by similar method, proceeding from the corresponding well-known statements in the theory of real Lie algebras (specifying Ado's theorem in these cases).
\enddemo

Let's review briefly the situation with triangular $D_2$-Lie algebras. Unlike any complex solvable Lie algebras, the famous fixed point theorem in the real case takes place only for triangular Lie algebras. For any solvable Lie algebras it is not true (even for one-dimensional Abelian Lie algebras). Due to the concept of the triangular Lie algebra we will consider one question, which has a little more general character---it concerns also the general theory of Lie algebras.

When studying solvable complex groups and Lie algebras, an important role is played by the concept of Borel subgroup and the corresponding Lie subalgebra---the maximum connected solvable Lie subgroup and the Lie subalgebra. It is connected, in particular, with the fact that for complex solvable groups and Lie algebras the fixed point theorem takes place. In  the real case for any solvable groups and Lie algebras such statement, generally speaking, is incorrect. However it is true for triangular groups and Lie algebras.
There is the natural question--- how properties of triangular Lie algebras and any solvable Lie algebras in the complex case are connected. For example, whether it is true that, if $R$ is some solvable complex Lie algebras, and $T$---its maximal triangular subgroup, then the complexification of this Lie subalgebra equals to $R$. It turns out that it is incorrect and therefore it turns out, that in case of complex solvable Lie algebras we have stronger statements that when we consider them as real Lie algebras.

Here we give some elementary example. Let's consider the complex Lie algebra of the type $K = \bold C + \bold C^n$ ---the semidirect sum, which is defined by a homomorphism $\phi: \bold C \to gl(n,\bold C) $. Let's consider matrix $A = \phi(1)$ and we will assume, that its eigenvalues are linearly independent over $\bold R$ (for example, when $n=2$, it means that the quotient of two eigenvalues of the corresponding matrix $A$ is not real number). Let $T$ be the  maximum triangular subgroup in such $R$. Let's prove that $T = \bold C^n$, what will means that the complexification of the Lie algebra of $T$ is not equal to $R$.

It is clear, that $\bold C^n \subset T$. Let $X \in K\setminus \bold C^n$. Let's prove that $X$ does not belongs to $T$. For this purpose we will consider eigenvalues of the operator $ad_X$ for its action on $\bold C^n$. It is clear, that they have to be proportional to the matrix $A$ eigenvalues. But then, as appears from our choice of matrix $A$, eigenvalues of $ad_X$ cannot be real. So it is proved that $X$ does not belongs to $T$ and therefore we have $T = \bold C^n$.

Similar questions arise also for $D_2$-Lie algebras. If $T$ is some  maximal triangular Lie subalgebra in the complex Lie algebra $R$, whether its complexification will equals to $R$ and whether it is possible to tell something similar about its dualization? Answers to both of these questions are similar to stated above--- the are negative. For this purpose it is enough to consider the Lie algebra $R \otimes D_2$, where $R$---the complex Lie algebra of the form $\bold C + \bold C^n$, constructed above. Here operator $\Cal E$ is nontrivial. There is a simpler example (though degenerated)---to take the Lie algebra $R$ with trivial $\Cal E$.

\head \S 2 About classification of approximate Lie algebras
\endhead

A classification of Lie algebras over $D_2$ is in many respects differs from a classification of real or complex Lie algebras. Let's begin with the fact, that there are no "purely" semi-simple Lie algebras here. Further, here can be several non isomorphic $D_2$-structures on one real Lie algebra, i.e. there can be non isomorphic $D_2$-Lie algebras, which realifications are isomorphic. It strongly expands the classification list of $D_2$-Lie algebras.

It is possible to carry out a classification, for example, by number $d(L)$ of essential parameters. It is obvious that for $d (L) $ =1 there are only two (up to $D_2$-isomorphism) $D_2$-Lie algebras, they are Abelian: one-dimensional $\bold R$ (with trivial $\Cal E$) and two-dimensional $D_2$ (isomorphic to $\bold R^2$ as the real Lie algebra).

There are known classifications of $D_2$-Lie algebras for cases of 2 and 3 essential parameters (see [3,4]). There are reason to believe, that actually $D_2$-Lie algebras in sense of our definition, at which the number of essential parameters is equal 2 or 3, are classified there (though concepts of $D_2$-Lie algebras are not used there). We will proceed from such understanding of articles [3,4].

It turns out that only since dimension 6 (and with three essential parameters) there are unsolvable $D_2$-Lie algebras. These are two such Lie algebras  $su(2) \otimes D_2$ and $sl (2, \bold R) \otimes D_2$ (both are of the type $S \otimes D_2 $, where $S$ --- one of two real simple three-dimensional Lie algebras). There are no classifications for bigger number of essential parameters and for bigger dimensions, because for this purpose, in particular, it is necessary to have classifications of real Lie algebras of dimensions 8 and more, that looks now as a unavailable task. The top of today's achievements in this direction is the classification of all Lie algebras of dimensions $\le 6$ and of nilpotent Lie algebras for dimensions up to 7. Therefore the result of article [3] is maximal possible now (it is possible to classify also all seven dimensional nilpotent $D_2$-Lie algebras, but the advantage of such work is not clear).

Let's provide the list of nilpotent $D_2$-Lie algebras of real dimensions $\le 3$ (considered up to $D_2$-isomorphism)---they are Abelian $\bold R, \bold R^2$ (with trivial $\Cal E)$, $D_2$ and non-Abelian $N (3, \bold R), D_2 \oplus \bold R$. For cases of dimensions of 4,5 and 6 see [3].

In classifications of $D_2$-Lie algebras, which can be useful for applications, it is useful to distinguish different classes of $D_2$-Lie algebras, because Lie algebras, which are of interest, belongs  sometimes only to one of this classes. Below are described some classes of $D_2$-Lie algebras (all of them are supposed real and finite-dimensional):

1. Any $D_2$-Lie algebras

2. $D_2$-Lie algebras of the form $L \otimes D_2$ for various "ordinary" Lie algebras $L$

3. The Lie subalgebras in Lie algebras of the type $L \otimes D_2$, which are $D_2$-modules

4. $D_2$-free Lie algebras  (i.e. $D_2$-Lie algebras, being free as $D_2$-modules)

5. Nondegenerate $D_2$-Lie algebras.

Let's specify those relations, which take place between the listed five classes of $D_2$-Lie algebras.

First, classes 1 and 3 are coincide. It follows from the Proposition 8 (Ado's theorem). It is clear also, that the class 1 is the most general, it contains all others. It is necessary to find out interrelations of classes 2, 4 and 5. It is clear, that the class 2 contains in classes 4 and 5. The classes 4 and 5 are coincide --- it follows from the analysis of their definitions. But these two classes do not coincide with class 2. In other words, not any nondegenerate $D_2$-Lie algebra can be presented in the form of the dualization of some Lie algebra. The matter is that dualization gives to Lie algebras very special structure --- they are presented in the form of the semidirect sum of the initial Lie algebra and some Abelian ideal (and both summands have identical dimensions). For any nondegenerate $D_2$-Lie algebras such decomposition not always exists.

In conclusion of article---several remarks about $D_2$-Lie algebras $L$ of vector fields on the $D_2$-line (which is not the ordinary two-dimensional real plane $\bold R^2$, but supplied standard $D_2$-structure).

For the line $D_2$ there is one dimensional invariant distribution, generated by the subspace of purely imaginary elements $\epsilon b$. It is integrable and therefore all $D_2$-Lie algebras of vector fields on $D_2$ are imprimitive (i.e. they keep some one-dimensional foliations). There are only four primitive Lie algebras of vector fields on the plane $\bold R^2$, their dimensions are $\le 8 $ (about S. Lie's classification in this case see, for example, [6]). But imprimitive Lie algebras of vector fields on $D_2$ (including those, which have a $D_2$-structure), may have arbitrarily large dimension.

Let's note that classical classification of Lie of vector fields on $\bold R^2$ uses the concept of similarity of Lie algebras of vector fields, and their similarity as $D_2$-subalgebras does not follow from this classification. It is necessary to carry out additional considerations, that greatly increases problems of classification. However there is one useful simplification --- such Lie algebras of vector fields on the plane are surely imprimitive.

Also let's note what full classification of Lie algebras of vector fields on the dual plane $D_2^2$ would demand to have classification of Lie algebras of vector fields on $\bold R^4$, that gives us now an  unattainable task. Classifications of approximate Lie algebras of dimensions $\le 6 $ (actually--- of $D_2$-Lie algebras) on the plane see in [4].

\enddocument